\parindent=0in
\magnification=\magstep1

\def\d{\delta}
\def\e{\epsilon}

\def\l{\lambda}
\def\m{\mu}
\def\n{\nu}
\def\o{\omega}

\def\z{\zeta}

\def\N{\bf N}

\def\R{{\bf R}}
\def\C{{\bf C}}
\def\P{{\bf P}}

\def\i{\infty}
\def\I{\int}

\def\s{\sum}

\def\sq{{\sqrt{-1}\over 2\pi}}

\def\sub{\subseteq}
\def\ra{\rightarrow}

\def\G{\Gamma}

\def\cL{{\cal L}}
\def\cX{{\cal X}}
\def\cM{{\cal M}}
\def\cK{{\cal K}}

\def\v{\vskip .1in}

\def\[{{\bf [}}
\def\]{{\bf ]}}

\def\pl{\partial}
\def\ddb{\partial\bar\partial}

\def\Chow{{\rm Chow}}

\def\si{\sigma}

\v\v\v
\centerline{\bf STABILITY, ENERGY FUNCTIONALS,
AND}
\centerline{{\bf K\"AHLER-EINSTEIN METRICS }
\footnote*{Supported in part by
the National Science Foundation under grants
DMS-98-00783 and DMS-01-00410}}
\v\v

\centerline{D.H. Phong${}^{\dagger}$ and Jacob Sturm${}^{\ddagger}$}

\v\v

\centerline{${}^{\dagger}$ Department of Mathematics}
\centerline{Columbia University, New York, NY 10027}
\v
\centerline{${}^{\ddagger}$ Department of Mathematics}
\centerline{Rutgers University, Newark, NJ 07102}

\v\v\v
{\S 1.\bf Introduction}

\v\v

Starting with the work of Yau [Y1], Donaldson [D1],
and Uhlenbeck-Yau [UY],
the notion of stability has revealed itself under many
guises to be closely related to the existence of canonical
metrics
in K\"ahler geometry. The equivalence between
Hermitian-Einstein metrics
on vector bundles and Mumford
stability was proved by Donaldson and Uhlenbeck-Yau in
[D1] and [UY], while the existence of K\"ahler-Einstein metrics
was conjectured in the early 1980's by Yau [Y2] to be equivalent
to stability in geometric invariant theory. At the present time,
the Yau conjecture has been
at least partially confirmed. The existence of K\"ahler-Einstein
metrics
has been shown to imply $K$-stability and $CM$-stability
by Tian [T2], and more
recently to imply Chow-Mumford stability by Donaldson [D2].

\v

In moduli theory, a non-zero vector
$\Chow(A)$ in a vector space $\C^{N+1}$
is associated to geometric objects $A$ (such as
vector bundles or varieties). The vector $\Chow(A)$ is
defined up to multiplicative constants, and its
$GL(N+1)$ orbit inside $\P^N$ uniquely determines $A$.
Thus the moduli space can be
constructed as the space of orbits.
The vector $\Chow(A)$ is said to be stable
if $SL(N+1)\cdot \Chow(A) \sub \C^{\N+1}$ is closed and the
stabilizer of $\Chow(A)$ is finite.
Stability is of particular importance,
since Geometric Invariant Theory
guarantees that the space of stable orbits
has the structure of an algebraic
variety [Mu]. On the other hand,
the emergence of stability as a necessary condition
for the existence of canonical metrics can be
quite subtle. Basic to Tian's approach is a remarkable asymptotic
equivalence
between a Lagrangian for constant scalar curvature metrics,
namely the Mabuchi energy, and a norm $||\cdot ||_Q$
on the space of Chow vectors which he
constructed by $\pl\bar\pl$ methods
and identified with suitable Quillen
metrics [T1, T2]. Quillen metrics had been related to
Lagrangians for Hermitian-Einstein metrics by
Donaldson [D3]. The asymptotic equivalence 
between $||\cdot||_Q$ and the Mabuchi energy becomes
exact for hypersurfaces in $\P^{n+1}$.
Donaldson [D2] relies instead on the concept of balanced manifolds
and Lu's recent evaluation [Lu1] of the Tian-Yau-Zelditch
expansion for the Bergman kernel on positive line bundles.

\medskip

In this paper, building on the earlier work of
Yau [Y3, Y4], Tian [T1,T2]
and Zhang [Z],
we construct a semi-norm
$||\cdot||_{\#}$ which is defined on the
{\it full} space $H^0(Gr, O(d))$ of polynomials of degree $d$
on the Grassmann variety $Gr$, and which gives
{\it exactly} the Mabuchi functional when restricted to
the space of Chow vectors, up to a specific
current supported on
the singular locus of the Chow variety.
This new semi-norm can be described
quite explicitly: Let  $Gr=Gr(N-n-1,\P^N)$ be the
space of $N-n-1$ projective planes in $\P^N$,
$P\ell: Gr(N-n-1,\P^N)\rightarrow \P(\wedge^{N-n} \C^{N+1})$ be
the Pl\"ucker imbedding, and $O(1)=P\ell^*H$, where $H$ is the
hyperplane bundle on $\P(\wedge ^{N-n}\C^{N+1})$.
Let $\o_{Gr}=P\ell^*\o_{FS}$, where $\o_{FS}$ is the
Fubini-Study K\"ahler form on the space $\P(\Lambda^{N-n}\C^{N+1})$.
Let $d$ be an arbitrary positive integer, and denote by
$m+1=(N-n)(n+1)$ and $D=\int_{Gr}\o_{Gr}^{m+1}$
respectively the dimension and the volume of $Gr$.
Then for any $f\in H^0(Gr,O(d))$, we define the norm $||f||_{\#}$
of $f$ by
$$
\eqalign{
{\rm log}\,||f||_\#^2 \ =& \ {(m+1)\over (m+2)(d-1)}\cdot{1\over D}\I_Z
\log\left({\o_{Gr}^m \wedge \ddb{|f(z)|^2\over |P\ell (z)|^{2d} }
\over
\o_{Gr}^{m+1} } \right) \o_{Gr}^m
\cr
&\ \ +\ \  {d-m-2\over (m+2)(d-1)}\cdot {1\over D} \I_{Gr}
\log{|f(z)|^2\over |P\ell (z)|^{2d} } \ \o_{Gr}^{m+1}
\cr}
\eqno(1.1)
$$
where $Z = \{z\in Gr: f(z)=0\}$.
\v
One easily sees that $||\cdot||_{\#}$ defines a semi-norm on the
finite dimensional
vector space $H^0(Gr,O(d))$: That is,
$||\l\cdot f||_{\#} \ = |\l|\cdot ||f||_{\#} $ for every $\l \in \C$
and $f \in H^0(Gr,O(d))$. Moreover, $||f||_{\#} \geq 0$ for all $f$.
On the other hand, $||\cdot||_{\#}$ is not a norm since for
$d > 1$ there exist
non-zero elements $f \in H^0(Gr,O(d))$ such that that $||f||_{\#}=0$.
\v

Our main result can be described as follows.
Let $X \sub \P^N$ be a smooth algebraic variety of dimension $n$,
$Z \sub Gr$ be the corresponding Chow variety, and let
$\Chow(X)=f \in H^0(Gr, O(d))$ be a defining section for $Z$.
Note that $Z$ is a singular variety. We assume that
the embedding $X\sub \P^N$ is {\it generic}
in the sense defined in \S 6.
\v
Let $Z_s$ be the singular locus of $Z$, and let
$Y_s = \{(x,z)\in X\times Z_s;x\in z\}
\subset \P^{N}\times Gr$.
Let $\o_Z$ be the restriction of $\o_{Gr}$ to
$Z_0=Z\setminus Z_s$. Then
the Ricci curvature $Ric(\o_Z)$
of $\o_Z$ is a smooth
$(1,1)$ form on $Z_0$.
Let $s(\o_Z)$ be the scalar curvature of
$\o_Z$, $V = vol(Z)$,
and define
$$ \m(Z) \ = \ {1\over V}\I_{Z_0} s(\o_Z)\,\o_Z^m\
\eqno(1.2)
$$
Let $[Y_s]$ be the current corresponding to $Y_s$ (see \S
6 for the precise definition) and define
$$
\deg(Y_s) \ = {1\over V}\,\langle [Y_s],\o_Z^{m-1}\rangle
$$

\v

Let $\o$ be the restriction of the Fubini-Study metric to
$X$, and let $\nu_{\o}(\phi)$
be the Mabuchi energy on
$X$ (see
\S 5 for the precise definition). For each $\sigma\in GL(N+1)$,
let $\phi_{\sigma}$ and
$\Phi_{\sigma}$ be the following functions on $\P^N$
and $Gr$ respectively
$$
\phi_{\sigma}(z)={\rm log}\,{|\sigma z|^2\over |z|^2},
\quad\quad
\Phi_{\sigma}(z)={\rm log}\,{|P\ell(\sigma z)|^2\over |P\ell(z)|^2}
$$
Define a generalized Mabuchi energy
$\nu_{\o}^{\#}(\phi_\si)$ by
$$
\nu_{\o}^{\#}(\phi_\si)
=
\nu_{\o}(\phi_{\sigma})
+{1\over V}
\langle[Y_s],\Phi_{\sigma}\sum_{i=0}^{m-1}\o_Z^i\sigma^*\o_Z^{m-1-i}\rangle
-{D\over V}
\cdot{m\, {\rm deg}(Y_s)\over m+1}
\cdot {\rm log}\,{||\sigma\cdot \Chow(X)||^2
\over ||\Chow (X)||^2}
$$
where $||\cdot||$ is the norm defined in (4.1) below. Then

\v\v

{\bf Theorem 1.} {\it For $\sigma \in SL(N+1,\C)$ we have
$$
\nu_{\o}^{\#}
(\phi_\si)\
= \
{D(m+2)(d-1)\over V(m+1)}\,{\rm log}\,{||\sigma\cdot \Chow(X)||_{\#}^2
\over ||\Chow(X)||_{\#}^2}
\eqno(1.3)
$$
}

\v\v
An interesting new notion emerges from the proof of Theorem 1,
namely the generalized Mabuchi energy
$\nu_{\o_Z}^{\#}(\Phi_{\si})$
of the {\it singular} Chow variety
$Z$ (see (6.3) for its precise definition).
The Chow variety $Z$ contains a
singular locus $Z_s$ of codimension $1$, and the
generalized Mabuchi
energy $\nu_{\o_Z}^{\#}(\Phi_{\si})$
is defined accordingly as consisting
of the usual Mabuchi energy $\nu_{\o_Z}(\Phi_{\si})$
of the regular
part $Z_0=Z\setminus Z_s$,
together with additional current terms due to
$Z_s$ (or more precisely, $[Y_s]$).
Theorem 1 is then the consequence of two results,
which may be interesting
in their own right
(c.f. Lemmas 6.1 and 6.2). The first
result is that the right hand side of (1.3)
can be equated with
the generalized Mabuchi energy
$\nu_{\o_Z}^{\#}(\Phi_{\si})$
of the {\it singular} Chow variety $Z$.
The second result is that the Mabuchi energy $\nu_{\o}(\phi_{\sigma})$
of the projective variety $X$ can be equated with the Mabuchi
energy $\nu_{\o_Z}(\Phi_{\si})$ of the {\it regular} part $Z_0$
of its Chow variety $Z$.

\v\v

As mentioned above, a similar formula to (1.3)
with an asymptotic bounded
error has been proved by Tian [T1], with $||\cdot||_{\#}$ replaced by
the Quillen metric $||\cdot||_Q$ of a certain virtual line
bundle. In the case where $X\sub \P^N$ is a smooth
hypersurface, the
version of (1.3) with $||\cdot||_Q$ also becomes exact, and it
is likely that the two notions $||\cdot||_{\#}$ and $||\cdot||_Q$
coincide. This issue is more complicated for higher codimensions,
not just because the error term in the version of (1.3) with $||\cdot||_Q$,
but also because $||\cdot||_Q$ is presently defined only on the
space of Chow vectors, and not yet on the whole
of $H^0(Gr,O(d))$.

\v\v

Besides the fact that it is completely explicit
and satisfies the relation (1.3) exactly,
the semi-norm $||f||_{\#}$ in (1.1) has several attractive features
which may be valuable in future investigations of
the relation between various notions
of stability and the
existence of K\"ahler-Einstein metrics. Indeed,
the existence of K\"ahler-Einstein metrics is known to imply
the boundedness from below of energy functionals
(see Siu-Yau [SY], Bando-Mabuchi [BM], and Ding-Tian [DT]).
Theorem 1 suggests that the boundedness from below of the Mabuchi
energy functional can eventually be related to the boundedness from below
of the norm $||f||_{\#}$, which is
defined even for singular
varieties. Furthermore, $||f||_{\#}$ is degenerate, so that
the condition that $||\sigma\cdot \Chow(X)||_{\#}/||\Chow(X)||_{\#}
\to \infty$ appears to be a stronger notion of stability than
the usual notion of Chow-Mumford stability.
It may be closely related to the notion of $CM$-stability
introduced earlier by Tian [T1].

\v\v

Our approach is based on an exact evaluation of the derivative
of  ${\rm log}\,||\sigma(t)\cdot \Chow(X)||$ along each
$1-parameter$ orbit
of $SL(N+1)$. This method appears to be technically
simpler than the approach in [T1][T2], which is based instead
on the evaluation of $\pl\bar\pl\,{\rm log}\,||\sigma\cdot \Chow(X)||$
on $SL(N+1)$.
\v
Our method applies equally well to
other contexts, namely to the component $L(h,k)$ of
the Donaldson functional ${\bf L}(k,h)$ for the
existence of Hermitian-Einstein metrics on vector bundles,
and to the component $F_{\o}^0(\phi)$ of the Lagrangian
$F_{\o}(\phi)$ for K\"ahler-Einstein metrics. In the case
of $F_{\o}^0(\phi)$, an exact relation of the form (1.3) had been
obtained by Zhang [Z], using the theory of Deligne pairings
[D]. Asymptotic
relations modulo $O(1)$ terms had been obtained by Paul [P]
and Wang [W] respectively for $F_{\o}^0(\phi)$ and for
$L(h,k)$. We shall use our approach to
give a unified and simpler proof of these earlier
results. We present these results in
(\S 2-\S 4) before proceeding to the more complicated calculation
of the Mabuchi energy, and take the opportunity to
mildly strengthen  earlier results of Zhang [Z], Luo [Luo],
and Wang [W] whenever it is readily
possible to do so by our methods.
In particular, we eliminate the $O(1)$ error terms in [W]
by showing that $L(h,k)$, restricted to
$G = SL(N, \C)$, equals exactly the log of the Gieseker point
in a suitable norm.
Wang [W] had shown that a vector bundle can be uniquely balanced if and only if its Gieseker point is stable.
Similarly, we show
that a manifold can be uniquely balanced if and only if
its Chow point
is stable.

\v\v\v
{\S 2.  \bf The  Donaldson energy functional and balanced bundles}

\v\v

Let $(X,\o)$ be a K\"ahler manifold and $\pi: E \ra X$
a vector bundle of rank $r$. The Donaldson functional
${\bf L}(h,k)$ is defined as follows. Let $h,k$ be two
hermitian metrics on $E$. Since the space of hermitian metrics
is convex, we can connect $k$ to $h$ be a smooth
path of hermitian metrics $h_t$, $0\leq t\leq 1$, $h_0=h$,
$h_1=k$. Let $R_t$ be the curvature of $h_t$. Then the
Donaldson functional ${\bf L}(h,k)$
is defined by
$$
{\bf L}(h,k)
=
\int_0^1dt\int_X i\,tr(h_t^{-1}\pl_th_t\cdot R_t){\o^{n-1}\over
(n-1)!}
-{c\over V}\int_X{\rm log}\,({\rm det}\,(k^{-1}h)){\o^n\over n!}
$$
where the constant $c$ is given by
$c={2\pi n\over r}\int_X c_1(E)\wedge\o^{n-1}$.
In this section,
$V=vol(X)=\int_X\o^n$ denotes
the volume of $X$ with respect to the
K\"ahler form $\o$.

\v\v
Assume  $E^*$ is generated by sections
$s_1,..., s_N \in H^0(X, E^*)$. If $e \in E$ then
$s(e) = (s_1(e),...,s_N(e) ) \in \C^N$ so
$s\times \pi: E \hookrightarrow \C^N \times X$.
Let $h$ be the metric
on $E$ defined by $h(e) = \s |s_i(e)|^2 = |s(e)|^2$,
where $|\cdot |$ is the usual norm on $\C^N$.
\v
Fix $\{\tau_1,...,\tau_k\}$, a basis for $H^0(X, det(E^*))$.
Let $T$ be the
matrix $T = (a^\m_{i_1\cdots i_r}) \in M $,
 defined by the equation

$$   s_{i_1} \wedge \cdots \wedge s_{i_r} \ = \
\ s_{i_1 \cdots  i_r} \ = \
\s_\m\  a^\m_{i_1\cdots i_r} \tau_\m
$$
where $M$ is the space of matrices $(A^\m_{i_1\cdots i_r})$
with
$1 \leq i_1 < \cdots i_r \leq N$ and $1 \leq \m \leq k$.
The Gieseker point of $s = (s_1,..., s_N)$ is the
point $[T] = (a^\m_{i_1\cdots i_r}) \in \P (M) $.
Then $[T]$ uniquely determines the image of $E$ in $\C^N\times X$.
\v
For $\sigma \in G$ we shall write $s^\sigma = (s_1,...,s_N)\sigma$,
and $h^\sigma(e) = |s^\sigma(e)|$. We define

$$  L(\sigma) \ = \ L(h,h_\sigma) \ =  \
{c\over V}\I_X \log\left( {det(h^\sigma) \over det(h) } \right)\ {\o^n\over n!}
$$
which is the second term in the definition of the Donaldson
functional ${\bf L}(k,h)$.
In local coordinates, if $e_1,..., e_r$ is a basis of smooth
sections of $E\sub \C^N \times U$, then for each $m$ such that
 $1\leq m \leq r$, and
for each $x \in U$, $e_m(x) = e_{im} (x)$
is a column vector in $\C^N$. Hence $A(x) = (e_{im}(x))$
is an $N \times r$ matrix, and
$h$ is the  $r\times r$ matrix $h =\ ^t\bar A A$. Thus
$det(h) = \s_{i_1< \cdots < i_r} |det (A_{i_1 \cdots  i_r})|^2$.
Since $s_i(x)  = \s e_{im}(x) e_m^*$ we get
$det(h) = \s_{i_1< \cdots < i_r} |s_{i_1 \cdots  i_r}|^2$.
\v
Thus
$$  L(\sigma) \ = \ {c\over V}\I_X \log\left(
{\det(A(x)^*\sigma ^*\sigma A(x)) \over
\det(A(x)^* A(x)) }\right)\ {\o^n\over n!}
\ = \
{c\over V}
\I_X \log\left( {\s |
(a^\m_{i_1\cdots i_r})^\sigma \tau_\m (x)|^2 \over
\s |a^\m_{i_1\cdots i_r} \tau_\m (x)|^2
 } \right)\ {\o^n\over n!}
\eqno(2.1)
$$
where $A^*$ is the conjugate transpose of $A$,
$(a^\m_{i_1\cdots i_r})^\sigma$
is the natural action of $G$ on $M$, and the
summation is over all $\m$ with $ 1 \leq \m \leq k$ and all
$1 \leq i_1 < \cdots < i_r \leq N$.
We note that the integral in (2.1) is finite.
This suggests defining the following norm on the
vector space $M$: Set for each $a=(a^\m_{i_1\cdots i_r}) \in M$
$$
{\rm log}\,||a||^2
=
{1\over V}\int_X \log \s |a^\m_{i_1\cdots i_r} \tau_\m (x)|^2\
{\o^n\over n!}
$$

\v
It is easy to see that $||\cdot||$ is a continuous norm, and hence
bounded on any compact subset of $M$.
In terms of $||\cdot||$, the formula (2.1) can be restated as

\v\v

{\bf Theorem 2.} {\it
Let $E \sub \C^N\times X$ be a vector bundle of rank $r$,
let $T$ be the Gieseker point of $E$, and
let $h$ be the metric on  $E$ defined by $h(e) = |e|_{\C^N}$.
For $\sigma \in G = SL(N, \C)$, let $h_\sigma$ be the metric
on $E$ defined by $h_\sigma(e) = |\sigma(e)|_{\C^N}$. Then
$$
L(h,h_\sigma) \ = c\,\log {||\sigma T||^2\over ||T||^2}
$$

}

\v\v

{\it Remark.} This is slightly more precise
than a result of Wang [W], who shows that
for any norm $||\cdot||$ on $M$, one has
$
L(h,h_\sigma) \ \geq  \ c\, \log\left( ||\sigma T||^2 \right) + C
$
for some constant $C$.
\v

\v
According to the theorem of Kempf-Ness [KN],  $E$ is
Gieseker stable if and only if $||\sigma T||^2$
is a proper map from $G$ to $\R$ (the inverse image
of a compact set is compact). This is equivalent to
requiring that $\log||\sigma T||^2$ is bounded below
by a positive constant and that
$\lim_{\sigma \to \i} \log||\sigma T||^2 = \i$
(ie, for every $B>0$ there exists a
compact subset $K \sub G$ such that $L(\sigma) \geq B$
if $\sigma \notin K$).
Combining
this with Theorem 2 we have the following corollary:

\v\v

{\bf Corollary 2.1} {\it Let $E \sub \C^N\times X$ be a vector
bundle of rank $r$. Then $E$ is Gieseker stable if
and only if the following conditions hold:

1. $L(\sigma) \geq \e > 0$ for some $\e$.

2. $\lim_{\sigma \to \i} L(\sigma) = \i\ .$
}

\v\v

{\bf Definition 2.1} {\it We say $E \sub \C^N\times X$ is balanced
if
$$  {1\over V}\I_X A(x) A(x)^* \ dV \ = \
{r \over N} \cdot I_{N\times N}
$$
where $A(x) = (a_1, \cdots a_r)$ is an orthonormal basis
of $E_x \sub \C^N$. We say $E$ can be (uniquely) balanced if and
only if there exists a (unique) $\sigma_0 \in SU(N)\backslash G$
such that $\sigma_0(E)$ is balanced.}

\v\v

{\it Example.} Let $X = Gr(r,N)$, the Grassmannian variety of all
$r$ planes in $\C^N$. Let $E$ be the canonical vector bundle
on $X$ of rank $r$. Then one easily sees that $X$ is balanced.

\v\v
{\bf Theorem 3.} (Wang) {\it Let
$E \sub \C^N\times X$ be a vector bundle.
Then the bundle $E$ can be uniquely balanced if
and only if its Gieseker point is stable.
}
\v\v

Theorem 3 as well as Lemmas 2.1 and 2.2 below
are due to Wang [W],
under a slightly different
formulation. He uses arguments from
the theory of moment maps. Here we shall
provide a direct calculation, along the lines
followed later for the proof of Theorem 1.

\v\v
{\bf Lemma 2.1 }   {\it $L$ has a critical point
at $\sigma_0 \in G$ if and only if $\sigma_0(E)$ is balanced.}
\v
{\it Proof.}
Write
$\sigma(t) = exp(tc)\sigma_0$ where $\sigma_0 \in G$ is
fixed, and $c$ is traceless.
Then
$$ {d \over dt}\log\det(A^*\sigma^*\sigma A) \ = \
tr((A^*\sigma^*\sigma A)^{-1}(A^*\sigma^*(c^*+c)\sigma A)
\eqno(2.2)
$$
Now $\sigma_0$ is a critical point if and only
if ${d \over dt}L(exp(ct)\sigma_0) = 0$ for all traceless $c$.
Replacing $A$ by $A\sigma_0$ we may assume that $\sigma_0=I$.
We may also choose our local sections $e_1,...,e_r$
to be orthonormal. Then $A^*A = I$ and (2.2) implies that
$\sigma_0$ is critical if and only if
$$
\I_X tr \big(A^*uA\big) \ dV \ = \ 0
$$
for all traceless hermitian $u$. This condition
is equivalent to $E$ being balanced. Q.E.D.
\v

Next we  show  $L(\sigma)$ is convex
in the following sense: View $L(\sigma)$ as a function on
the symmetric space
$SU(N)\backslash G$. Then $SU(N)\backslash G$ is a Riemannian manifold
whose geodesics are of the form $ exp(tc)\sigma_0$ where
$c$ is a traceless matrix with the property $c^* = c$.

\v\v

{\bf Lemma 2.2 } {\it Let $\sigma_0 \in G$ and let $c$ be
an arbitrary traceless self-adjoint matrix. Then
$$
{d^2\over dt^2}L(exp(tc)\sigma_0) \geq 0
\eqno (2.3)
$$}

{\it Proof.} Differentiating (2.2) again gives
$${d^2 \over dt^2}\log\det(A^*\sigma^*\sigma A) \ = \
tr\Big((A^*\sigma^*\sigma A)^{-1}
\big(A^*\sigma^*[c^*(c^*+c)+(c^*+c)c]\sigma A\big)\Big)
$$
$$
\ - \ tr\Big( (A^*\sigma^*\sigma A)^{-1}
(A^*\sigma^*(c^*+c)\sigma A)
(A^*\sigma^*\sigma A)^{-1}
(A^*\sigma^*(c^*+c)\sigma A)\Big)
$$
We may assume that
$\sigma_0 = I$ and $A^*A = I$. Let $u = c + c^* = 2c$.
Thus $u = u^*$ and $u$ is traceless.
The preceding equation becomes
$$
{d^2 \over dt^2}\log\det(A^*\sigma^*\sigma A) \ = \
tr\Big(
A^*[u^2] A\Big)
\ - \ tr\Big(
(A^*(u) A)
(A^*(u) A)\Big) \ = \
$$
$$
tr(A^*u(1-AA^*)uA) \ = \
tr\big((1-AA^*)uAA^*u^*\big)
$$
This last trace is non-negative since $uAA^*u^* \geq 0$ and
$1-AA^* \geq 0$. To see this last inequality, let $\l$ be
an eigenvalue of $AA^*$. Then $AA^*v = \l v$ for some
non-zero vector $v$. Applying $A^*$ to both sides, and
using the fact that $A^*A = I$, we get $A^*v = \l A^* v$.
Thus $\l = 1$ or $\l = 0$. Hence all the eigenvalues of
$1- AA^*$ are $\geq 0$ and $1-AA^* \geq 0$.

\v\v

{\it Proof of Theorem 3. }
Assume first that $L$ has a critical point $\sigma_0$. If
$\sigma \in G$ is any other point, then we can join
$\sigma_0$ to $\sigma$ by a geodesic. Lemma 2.2 implies
that $L$, restricted to the geodesic, has a minimum at
$\sigma_0$. Thus $L(\sigma) \geq L(\sigma_0)$, so $L$ is
bounded below. Theorem~2 implies that $||\sigma T||$ is
bounded below, and thus by definition the
Gieseker point is semi-stable. If the critical point $\sigma_0$ is unique,
then $L(\sigma)$ achieves its minimum at $\sigma_0$
and the same argument shows that $L$, restricted to
any geodesic through $\sigma_0$, goes to infinity. Thus
the Gieseker point is stable (by virtue of the one parameter
subgroup criterion)
\v
As for the converse, assume  that $E$ is stable.
Then the Kempf-Ness theorem says that $\log||\sigma T||$ is
proper on $H$ and bounded below. Thus it has a critical
point. The critical point is unique, for otherwise,
$\log ||\sigma T||$ would be constant on the
geodesic joining two critical points and thus it
could not be proper.

\v\v\v
{\S 3.\ \bf  The $F_{\o}^0$ functional and balanced manifolds}

\v\v
Consider now a smooth projective variety $X\subset\P^N$ of degree $d$ and dimension
$n$. Let $Z$ be the set of $(N-n-1)$-dimensional planes in $\P^N$ which
intersect $X$. Then $Z$ is contained in $Gr=Gr(N-n-1,\P^N)$ and has codimension
$1$. Thus there exists a holomorphic section $f\in H^0(Gr, O(d))$, unique
up to scalar multiplication,
which vanishes precisely on $Z$.
The section $f$ defines then a point in $\P H^0(Gr,O(d))$, which is
called the Chow point of $X\hookrightarrow \P^N$, and usually
denoted by $\Chow(X)$.

\v\v

We shall apply the same method as in the previous
section to prove the following:
\v\v
{\bf Theorem 4.} {\it  Let $X\sub \C P^N$ be a smooth projective
variety.
Then the Chow point of $X$ is stable if and only if
there is a unique $\sigma_0$
in  $\ SU(N+1)\backslash SL(N+1,\C)$
such that $\sigma_0(X)$ is balanced, i.e.,
$$ {1\over vol(X)}\I_{\sigma_0 (X)} \left(
{\bar z_jz_i \over |z_0|^2 + \cdots |z_N|^2}\right)
\o_{FS}^n \ = \ {1\over N+1}\cdot \d_{ij}
\eqno(3.1)
$$}

\v\v
{\it Remark.} This theorem is a mild strengthening
of a result of Zhang [Z],
who shows that balanced implies semi-stable and stable implies
uniquely balanced.
It also mildly
strenghthens a similar theorem of Luo [Luo].
\v\v

To prove Theorem 4, we make use of the component
$F_{\o}^0(\phi)$ of a Lagrangian $F_{\o}(\phi)$
for K\"ahler-Einstein metrics
$$
F_{\o}(\phi)
=F_{\o}^0(\phi)-{\rm log}\,({1\over V}
\int_Xe^{h_{\o}-\phi}\o^n)
$$
Here $\o$ is a K\"ahler metric on $X$,
$\phi$ is a smooth function in the K\"ahler cone of $X$,
$\o_{\phi}=\o+{\sqrt{-1}\over 2\pi}\pl\bar\pl\phi$,
$V=vol(X)=\int_X\o^n$,
and $h_{\o}$ is defined by $Ric(\o)-\o={\sqrt{-1}\over 2\pi}\pl\bar\pl h_{\o}$,
$\int_Xe^{h_{\o}}\o^n=\int_X\o^n$.
The component $F_{\o}^0(\phi)$
is defined by
$$
F_\o^0(\phi) \ = \ J_\o(\phi) - {1\over V}\I_X\phi\, \o^n
$$
where the functional $J$ is given by
$$
J_\o(\phi) \ = \ {\sqrt{-1}\over 2\pi V}\cdot {1\over n+1}\s_{i=0}^{n-1}
\ (i+1)
\ \pl \phi\wedge \bar \pl \phi \wedge\o_\phi^{n-i-1}\wedge \o^i
$$
The variational derivative of the functional $J_\o(\phi)$ is well-known
(c.f. [T3]).
For a smoothly
varying family of potentials $\phi(t)$, $t \in (-\e, \e)$, we have
$$
-{d \over dt}F_\o^0(\phi(t)) \ = \
{1\over V}\I_X\ \dot\phi(t)\ \o_\phi^n
\eqno(3.2)
$$
Henceforth we fix $\o=\o_{FS}
=\sq \pl\bar\pl\,{\rm log}\,|x|^2$ to be the Fubini-Study
K\"ahler form on $X$.
For
$\sigma \in G$ let
 $$\phi_\sigma  \ = \ \log \left({|\sigma x|^2 \over |x|^2}\right)
\ = \
\log \left({^t\bar x^t\bar \sigma \sigma x\over |x|^2}\right)
\ = \
\log \left({x^*\sigma^*\sigma x\over x^*x}\right) \
$$
Then
$$
\sigma^*\o =
\o+{\sqrt{-1}\over 2\pi}
\pl\bar\pl\,\phi_{\si}
\equiv
\o_{\phi_\sigma} \equiv \o_\sigma.
$$
Let $c$ be a traceless hermitian matrix, $\sigma(t) = exp(ct)\sigma_0$,
and $F(\sigma_0, c, t) = F_\o^0(\phi_{\sigma(t)})$. We say that
$\sigma_0$ is a critical point of $F_{\o}^0$ if and only if
$F'(\sigma_0, c, t)|_{t=0} = 0$
for all traceless hermitian $c$. Note that
$$
{d \over dt}\phi_{\sigma(t)} \ = \ \dot\phi_\sigma \ = \
{ { x^* \sigma^*(c^*+c)\sigma x\over
x^* \sigma^*\sigma x } }
\eqno(3.3)
$$

{\bf Lemma 3.1} {\it
\vskip .05in
1. A matrix $\sigma_0$ is a critical point of $F^0$ if and
only if the manifold $\sigma_0(X)$ is balanced.

2. For all $\sigma_0, c$  and $t$ we have the formula
$$ -F'(\sigma_0, c, t)\ = \
{1\over V}\I_X { x^* \sigma^*(c^*+c)\sigma x\over
x^* \sigma^*\sigma x } \ \o_\phi^n
\eqno(3.4)
$$
3. For all $\sigma_0, c$ and $t$ we
have $-F''(\sigma_0, c, t) \geq 0$ . In fact we have
$$-VF''(\sigma_0, c, t) \ = \
{1\over n!}
\I_X\iota_{\tilde c} (\omega_\phi^{n+1})
\eqno(3.5)
$$

where $\tilde c$ is the vector field on $\P^N$ generated
by the infinitesimal action of $exp(ct)$.}
\v
{\it Remark.}  The only non-trivial part of the lemma is
part 3.  Zhang [Z] also proves that $-F''(\sigma_0, c,t)\geq 0$,
by making use of  a result of Deligne [D].
Another proof was given by Tian [T4]. We give
a more elementary proof by a direct calculation.

\v
{\it Proof of Lemma 3.1} Part two follows immediately from
(3.2) and (3.3a), and part one follows from part two.
To prove part 3, we
differentiate both sides of (3.4):
$$
-VF'' \ = \
\I_X \Bigg\{ {(x^* \sigma^*\sigma x)
(  x^* \sigma^*[c^*(c^*+c)+ (c^*+c)c]\sigma x) -
(x^*\sigma^*(c^*+c)\sigma x)^2
\over (x^* \sigma^*\sigma x)^2} \o_\phi^n \
$$
$$ -  n \pl\left({ x^* \sigma^*(c^*+c)\sigma x\over
x^* \sigma^*\sigma x }\right) \wedge
\bar\pl\left({ x^* \sigma^*(c^*+c)\sigma x\over
x^* \sigma^*\sigma x } \right) \o_\phi^{n-1}\Bigg\}
\eqno(3.5)
$$
We may assume that $c$ is a traceless diagonal matrix
with real entries. As before, we write
$u = 2c$, and we let $u_0,...,u_{N}$ be the
diagonal entries of $U$. Since we can view
$\P^N $ as the set of elements in $\P^{N+1}$
whose first entry is zero, we may also assume
(after replacing $N$ by $N+1$)
that $u_0 = 0$.
\v
We may assume
that $\sigma = I$ so that $\o = \o_{FS}$.
Recall that
on the coordinate chart
$U_0 = \{(1,z_1,...,z_N) \sub \C\P^N\}$,
the Fubini-Study
metric is given by:
$$ \o \ = \ {dz_i\wedge d\bar z_i\over 1+|z|^2} \ - \
{\bar z_idz_i\wedge z_j d\bar z_j \over (1+|z|^2)^2}
$$
Thus if $A = \sum_{i=1}^N a_i\pl_{z_i}$
is a tangent vector
at the point $x = (1,z_1,...,z_N)$, then
$$
\o (x)(A,\bar A) \ = \
{x^*x(|a_1|^2+\cdots+ |a_N|^2) -
|a_1\bar z_1+\cdots + a_N\bar z_N|^2 \over
(x^*x)^2}
\eqno(3.6)
$$
At $\sigma=I$ and $\o=\o_{FS}$, the first term in the integrand
on the right side of (3.5) is then immediately seen to
coincide with $\omega(x)(A_0,\bar A_0)$, where $A_0$
is defined to be the vector $A_0=(u_1z_1,\cdots,u_Nz_N)$.
Similarly, an explicit calculation gives
$$
\partial({x^*ux\over x^*x})(A)
=
{1\over (x^*x)^2}
\bigg(x^*x\,(A|\bar A_0)-(A|\bar z)(z|\bar A_0)\bigg)
=\o(x)(A,\bar A_0)
$$
Let $A_1,\cdots,A_n$ be now $n$ arbitrary tangent vectors at the
point $x$. Then, by the definition of the wedge product,
$$
\eqalign{
(\partial({x^*ux\over x^*x})
&
\wedge
\bar\partial({x^*ux\over x^*x})
\wedge
\omega^{n-1}\big)(A_1,\cdots,A_n)\cr
&
=
\sum_{j,k=1}^n(-1)^{j+k}
\omega(A_j,\bar A_0)\omega(A_0,\bar A_k)
\omega^{n-1}(\{A_p\}_{p\not=j},\{\bar A_n\}_{q\not=k})
\cr}
$$
But wedge powers of a $(1,1)$-form are given by
$$
\omega^K(A_1,\cdots,A_K,\bar A_1,\cdots,\bar A_K)
=
K!
\
{\rm det}\,\omega(A_j,\bar A_k)
$$
Substituting this in the full integrand on the right hand side
of (3.5), we obtain
$$
n!
\bigg(
\omega(A_0,\bar A_0)
{\rm det}_{1\leq p,q\leq n}\,\omega(A_p,\bar A_q)
+
\sum_{j,k=1}^n(-1)^{j+k+1}
\omega(A_j,\bar A_0)\omega(A_0,\bar A_k)
{\rm det}_{p\not=j\atop q\not=k}
\,
\o(A_p,\bar A_q)\bigg)
$$
The expression between parentheses
can be recognized as the expansion along
the first row (or the first column) of the determinant
of the $(n+1)\times (n+1)$ matrix $\omega(A_p,\bar A_q)$,
$0\leq p,q\leq n$.
Since $A_0$ is readily recognized as the vector field $\tilde c$,
Part 3 of Lemma 3.1 follows. Q.E.D.

\v\v
{\it Proof of Theorem 4.} Once Theorem 5 is available,
Theorem 4 can be proved in exactly the same manner as Theorem 3. Q.E.D.

\v\v
{\bf \S 4. The $F_{\o}^0$ functional and the Chow point}
\v
\v
Fix positive integers $n < N$ and let $Gr(N-n-1,\P^N)$
 be the space of $N-n-1$ projective
planes in $\P^N$. Then $Gr(N-n-1,\P^N)= Gr(N-n,\C^{N+1})$, the
 set of
$N-n$ vector subspaces of $\C^{N+1}$. Note that $G = SL(N+1,\C)$
acts on $Gr$ in a natural way.
Recall the notation introduced in \S 1, namely
$P\ell: Gr(N-n-1,\P^N) \ra \P(\Lambda^{N-n}\C^{N+1})$
is the Pl\"ucker embedding, $O(1) = P\ell^*H$, where
$H$ is the
hyperplane bundle on $\P(\Lambda^{N-n}\C^{N+1})$,
and
$\o_{Gr} = Pl^*\Omega_{FS}$ where $\Omega_{FS}$ is
the Fubini-Study metric on $\P(\Lambda^{N-n}\C^{N+1})$.
\v
Now let $d$ be a positive integer and define a norm
on the vector space $H^0 (Gr, O(d)) $ as follows: If
$f \in H^0 (Gr, O(d)) $ then
$$
\log ||f||^2 \ = \ {1\over D} \I_{Gr}
\log{|f(z)|^2 \over |Pl(z)|^{2d} } \ \o_{Gr}^m
$$
where $D = \I_{Gr} \o_{Gr}^{m+1}$ and $m+1 = (N-n)(n+1)$ is the
dimension of $Gr$.

\v\v

{\bf Theorem 5.}
{\it Let $X \sub \P^N$ be a smooth projective variety
of dimension $n$. Let $\o$ be the Fubini-Study
metric on $X$ and let $f=\Chow(X)$ be the Chow point of
$X$. Then}
$$ -V(n+1) F^0_\o(\phi_\sigma) \ = \
 \ \log\,{||\sigma\cdot \Chow(X)||^2\over ||\Chow(X)||^2}
\eqno(4.1)
$$

\v\v
This theorem is equivalent to one proved by Zhang in [Z],
using Deligne pairings.
Paul [P] shows, using a different method, that the
difference of the left and right sides of (4.1) is a bounded
continuous function on $G$. He also shows how such estimates
can be applied to give a new proof of Mumford's theorem
on the stability of curves.
We shall give a different proof of
their theorem, in the spirit of the proof of Theorem 1.
\v
\v

{\it Proof.} Let
$\Gamma  = \{(x,z) \in \P^N \times Gr: x \in z\}$ and let
$\pi_i$ be the projection map of $\Gamma$ onto the two factors.
We make use of the following formula:
$$ {\pi_1}_*\pi^*_2 \o_{Gr}^{m+1} \ = \ D\o^{n+1}
\eqno(4.2)
$$
Since both sides are invariant under the $U(N+1)$ action,
they are equal up to a constant. The constant is the
ratio of $\I_{Gr}\o_{Gr}^{m+1}$ and $\I_{\P^{n+1}}\o^{n+1}$,
which is equal to $D$ by definition.
\v
Now let $\sigma(t)$ be a path in $G$, and set
$f^{\sigma}(z)=f(\sigma^{-1}(z))$. We compute

$$
\eqalign{{d \over dt}\log||f^\sigma||^2 \
=&
\
{1\over D}{d \over dt} \I_{Gr}
\log{|f(\sigma^{-1} z)|^2 \over |Pl(z)|^{2d} } \ \o_{Gr}^{m+1}
 = \
{1\over D}{d \over dt} \I_{Gr}
\log{|f( z)|^2 \over |Pl(\sigma z)|^{2d} } \ (\sigma^*\o_{Gr})^{m+1}
\cr
=&
\ {1\over D}{d \over dt} \I_{Gr}
\log{|f( z)|^2 \over |Pl( z)|^{2d} } \ (\sigma^*\o_{Gr})^{m+1}\ \ - \
{d\over D}{d \over dt} \I_{Gr}
\Phi_\sigma \ (\sigma^*\o_{Gr})^{m+1}\cr
\equiv&
\ A \ -\  B
\cr}
\eqno(4.3)
$$
where $\Phi_\sigma = \log\,{|Pl(\sigma(z)|^2/|Pl(z)|^2} $.
Since $\sigma^*\o_{Gr} = \o_{Gr} +\sq \pl\bar\pl \Phi_\sigma$, we
obtain, using the Poincare-Lelong formula
$$
\eqalign{A \ =
& \ {1\over D} \I_{Gr}
\log{|f( z)|^2 \over |Pl( z)|^{2d} } \ (m+1)(\sigma^*\o_{Gr})^{m}
\sq\pl\bar\pl \dot\Phi_\sigma\cr
=& \
{1\over D}\I_Z \dot\Phi_\sigma (m+1)(\sigma^*\o_{Gr})^{m} \
 - \
{d\over D}\I_{Gr} \dot\Phi_\sigma (m+1)(\sigma^*\o_{Gr})^{m} \o_{Gr}
\cr
B \ =&
\ {d\over D} \I_{Gr}
\dot\Phi_\sigma \ (\sigma^*\o_{Gr})^{m+1}\ \ + \
{d\over D}\I_{Gr}\Phi_\sigma\,\sq \pl\bar\pl \dot\Phi_\sigma
(m+1)(\sigma^*\o_{Gr})^{m}\cr}
$$
We integrate by parts in the second term in $B$.
Combining the result with the other term in $B$
and the second term in $A$, we can rewrite the
above equation as
$$
{d \over dt}\log||f^\sigma||^2 \ = \
{1\over D}\I_Z \dot\Phi_\sigma (m+1)(\sigma^*\o_{Gr})^{m} \ \ - \
{d(m+2)\over D}\I_{Gr}
\dot\Phi_\sigma \ (\sigma^*\o_{Gr})^{m+1}\
\eqno(4.4)
$$
After making a change a variables $z \mapsto \sigma^{-1}z$, we may assume
that the $\sigma$ in the second integral is
the identity.
We conclude that the second integral is zero since the hyperplane
bundle over the Grassmannian variety is balanced.
More simply, we may write $\dot\Phi_{\sigma}(d)=tr(Z^*\,d\,Z)$ with $Z^*Z=1$ and $d=c+c^*$, $tr(d)=0$.
Then the integral
$$
M(d)=\int_{Gr}tr(Z^*\,d\,Z)\o_{Gr}^{m+1}
$$
satisfies $M(u^*du)=M(d)$, $M(d_1+d_2)=M(d_1)+M(d_2)$.
The first property allows us to assume that $d$ is diagonal.
The second implies that $M(d)$ is the same if we average it
over the permutations of the eigenvalues of $d$. Thus we must
have
$$
M(d)=0.
$$

\v

It remains to show that the integral over $Z$ in (4.4)
can be reexpressed as an integral over $X$. For this,
we apply $\sigma(t) = exp(ct)\sigma_0$ to
both sides of (4.2),
differentiate with respect to $t$ and evaluate at $t=0$:
$$
{\pi_1}_*\pi^*_2
\big((m+1)\sigma^*\o_{Gr}^{m}\pl\bar\pl\dot\Phi_{\sigma_0,c}\big) \ = \
D(n+1)(\sigma^*\o^{n})\pl\bar\pl \dot\phi_{\sigma_0,c}
\eqno(4.5)
$$
where $ \Phi_{\sigma_0,c}$ ({\it resp.} $\phi_{\sigma_0,c}$)
is the derivative of $\Phi_{\sigma(t)}$
({\it resp.} $\phi_{\sigma(t)}$)
at
$t=0$.
This shows that
$$   {\pi_1}_*\pi^*_2
\big(\dot\Phi_{\sigma_0,c} (m+1)(\sigma_0^*\o_{Gr}^{m})\big) \ = \
D(n+1)\dot\phi_{\sigma_0,c}( \sigma_0^*\o^{n})  \ + \ \eta(\sigma_0, h)
\eqno(4.6)
$$
for some closed smooth form $\eta(\sigma_0, h)$ where $h = c+c^*$ is
traceless hermitian.
\v
We claim that $\eta(\sigma_0, h)$ is exact.
Note the following properties:   $\eta(\sigma_0, u^*hu) \ = \eta(u\sigma_0, h)$
for all unitary
matrices~$u$. Also, $\eta(\sigma_0, h) = \sigma_0^*(\eta(I, h))$ where
$I$ is the identity matrix, and $\eta(\sigma_0, h)$ is
a linear function of $h$. This shows that we may
assume that $\sigma_0 = I$ and $h$ is diagonal with real eigenvalues
whose sum is zero. In fact, using the linearity property
 we may assume that $h$ is the matrix
whose diagonal entries are $(1,-1,0,...,0)$. Now let
$\iota: \P^n \hookrightarrow \P^N$ be the map
$(x_0,...,x_n) \mapsto (x_0,...,x_n,0,0,..., 0)$ and let $X'$ denote
the image of $\iota$. To show that $\eta(I,h)$ is exact, it
suffices to show that its integral over $X'$ is zero.
\v
To see this, first observe that $\I_{X'} \dot\phi_{I,h}\o^{n} = tr(hM)$
where $M = \I_{X'} xx^* \o^n $ and $x$ is the vector of homogeneous
coordinates, normalized to have length one (this follows form
(1.2)). Let $U_n \sub U(N+1)$ be
the group of unitary matrices which map $X'$ into itself. Then $M$
commutes with $U_n$ which means that it's a diagonal matrix whose
first $n+1$ entries are all equal. But this means that $tr(hM)=0$
and thus $\I_{X'} \dot\phi_{I,h}\o^{n}=0$. Similarly,
$\I_{Z'}\dot\Phi_{I,h} \o_{Gr}^{m}  = 0$ where $Z'\sub Gr $
is the zero locus of the Chow point of $X'$. This
implies that  $\I_{X'} \eta(I,h)=0-0=0$,
which proves our claim.

\v

Plugging this into (4.4) we obtain:

$$
{d \over dt}\log||f^\sigma||^2 \ = \ \I_X (n+1)\dot\phi_\sigma
(\sigma^*\o^{n})
$$
Comparing this with the derivative of
$F_{\o}^0$ given in (3.2), we obtain Theorem 5. Q.E.D.

\v
{\it Remark.} A similar argument shows that for
every $k$ such that $0 \leq k \leq n$, there
is a constant $D_k> 0$ such that
$$ {\pi_1}_*\pi^*_2 \o_{Gr}^{m-k} \ = \ D_k\cdot\o^{n-k}
$$
$$    {\pi_1}_*\pi^*_2
\big(\dot\Phi_{\sigma_0,c} (m+1-k)(\sigma_0^*\o_{Gr}^{m-k})\big) \ = \
D_k(n+1-k)\dot\phi_{\sigma_0,c}( \sigma_0^*\o^{n-k})  \ +
\ \eta_k(\sigma_0, h)
\eqno(4.7)
$$
where
$\eta_k(\sigma_0, h)
$
is a smooth exact $(n-k,n-k)$ form.

\v\v

{\bf \S 5. The Mabuchi energy: Hypersurface Case}

\v\v

We turn now to the setting of main interest in the present
paper, namely that of Theorem 1.
In this section, we discuss first the case of hypersurfaces.
In this case, our construction of $||\cdot||_{\#}$ shares
many features with Tian's construction of $||\cdot||_Q$
in [T1]. However, the explicit formulas in our approach
will facilitate the generalization to higher codimension.
We present our proof of the hypersurface case
in such a way as to apply verbatim
to the higher codimension case as much as possible, so that
the difficulties inherent to this latter case
will be more transparent in the next section.

\v\v

We first recall the definition of the Mabuchi energy [M].
Let $X$ be a K\"ahler manifold, with K\"ahler form $\o$.
Then the Mabuchi energy $\nu_{\o}(\phi)$ is defined
for all K\"ahler forms $\omega_{\phi}=
\o+{\sqrt{-1}\over 2\pi}\pl\bar\pl\phi$ by choosing a path
$\o_t=\o+{\sqrt{-1}\over 2\pi}
\pl\bar\pl\phi_t$, $0\leq t\leq 1$, $\o_0=\o$, $\o_1=\o$,
and setting
$$
\nu_{\o}(\phi)
=
\ -\I_0^1\I_X \dot\phi_t (s(\o_t)-\m) \ \o_t^n\wedge dt
\ = \
\ -\I_0^1\I_X \dot\phi_t (n{\sqrt{-1}\over 2 } Ric(\o_t)-\m\o_t)\o_t^{n-1} \
\wedge dt
$$
where $\o_t = R(h_t) = -\sq\pl\bar\pl \log h_t$
 is the curvature of $\o_t$, the function
$s(\o_t)$ is the scalar curvature of $\o_t$, and
$\m = n(n+2-d)$
is the average scalar curvature. The Ricci curvature
is defined by $Ric(\o) = -\sq\pl\bar\pl \log (\o^n)$.
Thus if we view $\o^n$ as a metric on the canonical
bundle $K = \Lambda^nT^*_X$, and then the Ricci curvature
is just the curvature of this metric.

\v
There is another more direct way of defining the Mabuchi
energy, pointed out by
Tian [T1] and Chen [C1], which does not require
an integral along paths:
$$
\nu_{\o}(\phi)
=
{1\over vol(X)}\int_X\bigg({\rm log}\,({\o_{\phi}^n\over\o^n}) \o_{\phi}^n
+h_{\o}(\o^n-\o_{\phi}^n)\bigg)
-{1\over n}(I_{\o}(\phi)-J_{\o}(\phi))
$$
where $\o_{\phi}=\o+{\sqrt{-1}\over 2\pi}\pl\bar\pl\phi$,
$vol(X)=\int_X\o^n$ is the volume of $X$,
$J_{\o}(\phi)$ is the functional introduced before,
and $I_{\o}(\phi)={1\over vol(X)}\int_X\phi(\o^n-\o_{\phi}^n)$ is
the other functional introduced by Yau [Y1] and Aubin. In [Z], Zhang used the
Deligne pairing $<{\cal L}_0,{\cal L}_1,\cdots,{\cal L}_n>_{(X/S)}$
with ${\cal L}_i=O(1)$ for all $0\leq i\leq n$
to obtain the functional $F_{\o}^0$. Here $\pi:X\rightarrow S$
is a flat projective morphism of integral schemes of
relative dimension $n$, so that each fiber $X_s$, $s\in S$
is a projective variety of dimension $n$, and ${\cal L}_i$
are line bundles over $X$. If we choose instead
${\cal L}_i=O(1)$ for $0\leq i<n$, ${\cal L}_n={\cal K}$,
then we obtain another expression for the Mabuchi functional
which does not require a path integration
$$
\nu_{\o}(\phi)=-E_{\o}(\phi)-\mu V\,F_{\o}^0(\phi)
$$
Here the functional $E(\phi)$ is defined by
$$
<O(1)\otimes {\cal O}(\phi), \cdots,
O(1)\otimes {\cal O}(\phi), {\cal K}\otimes{\cal O}({\rm log}\,
{\o_{\phi}^n\over\o^n})>_{X/S}
=<O(1),\cdots,O(1),{\cal K}>_{X/S}\otimes {\cal O}(E)
$$
For our purposes, it is most convenient to observe that $F_{\o}^0$
can be rewritten as
$$
F_{\o}^0(\phi)
=-{1\over vol(X)}{1\over n+1}\int_X\phi\sum_{i=0}^n\o^i\o_{\phi}^{n-i}
\eqno(5.1)
$$
and to recast the Mabuchi energy $\nu_{\o}(\phi)$
in the following form
$$
\n_\o(\phi) \ = \ {1\over vol(X)}
\int_X \Bigg\{\log\left({\o^n_\phi\over \o^n}\right)\o_\phi^n
-\phi\Big( Ric(\o)\s_{i=0}^{n-1}\o^i\o_\phi^{n-1-i}
\ - \
{\m(X)\over n+1}\s_{i=0}^n
\o^i\o_\phi^{n-i}\Big)\Bigg\}
\eqno(5.2)
$$
where $\m(Z)$, the average of the scalar curvature of
$\o_X$, is given by the following formula:
$$ \m(X) \ = \ m(m+2-d)
\eqno(5.3)
$$
Here we use $vol(X)$ to distinguish the volume of $X$
from the volume $V=vol(Z)$ of the Chow variety $Z$,
which also enters our formulas.
\v\v
Consider now the case of a smooth hypersurface $X\subset\P^{n+1}$.
In this case, $X$ coincides with its Chow variety $Z$, $n=m$,
the Grassmannian $Gr$ reduces to $\P^{m+1}=\P^{n+1}$,
$vol(X)=vol(Z)=V$,
and the functions $\Phi_{\si}$ and $\phi_{\si}$ coincide.
Theorem 1 simplifies considerably,
and we restate it as follows

\v\v

{\bf Theorem 6.} {\it Let $Z\sub \P^{m+1}$ be a  smooth
hypersurface of degree $d$ and let $f$ be the
section of $ H^0(\P^{m+1}, O(d))$ which defines $Z$.
Let $\o_{FS}$ be the Fubini-Study metric on $\P^{m+1}$,
and let $\o$ be the restriction of $\o_{FS}$ to $Z$.
For $\sigma \in
SL(m+2,\C)$ we have
$$
\n_{\o}(\Phi_\sigma)\
= \
{D(m+2)(d-1)\over V(m+1)}\,{\rm log}\,{||\sigma\cdot f||_{\#}^2
\over ||f||_{\#}^2}
\eqno(5.4)
$$
}

\v\v

{\it Proof.}
We evaluate first the contributions
of the logarithmic terms
from $||f||_\#$ on the right side of (5.4)
$$
{1\over V}\cdot\I_Z \log
\left(
{\si^* \o^m \wedge \ddb {
|f(z)|^2 \over |\si z|^2
} \over
\si^*\o^{m+1}
}
\right)
\si^*\o^{m}
\ - \
{1\over V}\cdot\I_Z \log
\left(
{\o^m \wedge \ddb {
|f(z)|^2 \over | z|^2
} \over
\o^{m+1}
}
\right)
\o^{m}
\eqno(5.5)
$$
Writing $\o^m = (\o^m - \si^*\o^m) + \si^*\o^m$ in the second
integral, the right side of  (5.5) becomes:
$$
\eqalign{
&\I_Z \log
\left(
{\si^* \o^m \wedge \ddb {
|f(z)|^2 \over |\si z|^2
} \over
 \o^m \wedge \ddb {
|f(z)|^2 \over | z|^2
}
}
\cdot
{\o^{m+1}\over \si^*\o^{m+1}
}
\right)
\si^*\o^m
\cr
&
\quad\quad
-
\I_Z \log
\left(
{\o^m \wedge \ddb {
|f(z)|^2 \over | z|^2
} \over
\o^{m+1}
}
\right)
(-{\sqrt{-1}\over 2\pi }\ddb \Phi_\si)
\s_{i=0}^{m-1} \o^{m-i}\si^*\o^i
\equiv A-B\cr}
$$
To evaluate the integrands in the last equation,
we choose $(m+1)$ linearly independent
holomorphic vector
fields $X_1,...X_{m}, Y$ in an open neighborhood of
a point $z \in X \sub \P^{m+1}$ in such a way
that $X_1,...,X_m$ are tangent to $Z$. Using
the definition of the wedge product, we
evaluate
$$
\left(\sigma^*\o^m\wedge \ddb
\left({|f(z)|^2\over |\sigma z|^{2d} }\right) \right)
(X_1,X_2,...,X_m,Y) \ =
$$
$$
\ \s_\pi
\sigma^*(\o) \otimes \cdots \sigma^*(\o)\otimes \ddb
\left({|f(z)|^2\over |\sigma z|^{2d} }\right)
 (X_1\otimes\bar X_1\otimes X_2\otimes\bar
X_2\cdots X_m\otimes\bar X_m,Y\otimes\bar Y)^\pi
$$
where $\pi$ ranges over all the permutations $\pi$ of the sequence
$(X_1,\bar X_1,...,X_m,\bar X_m,Y,\bar Y)$. Since $f$ vanishes
on $Z$, the only terms in the sum which are non-zero are those
corresponding to  permutations which permute the $X_i$ and the
$\bar X_i$, but fix $Y$ and $\bar Y$. Furthermore,
note that when we apply $\ddb$ to
$\left({|f(z)|^2\over |\sigma z|^{2d} }\right)$, the product
rule will yield several terms. But when we restrict to $X$,
the only term which doesn't vanish is the one where
the $\ddb$ lands on $|f(z)|^2$ (again, since $f$ vanishes
on $Z$). Thus
$$
\ddb\left({|f(z)|^2\over |\sigma z|^{2d} }\right) (Y\otimes \bar Y)
\Bigg|_Z \ = \
{1\over  |\sigma z|^{2d} }\cdot|Y(f)|^2
\eqno(5.6)
$$
and we may write
$$
\eqalign{
A&=\I_Z \log
\left(
{\si^* \o^m(X_1,..., X_m) \cdot  {
|Y(f)|^2 \over |\si z|^{2d}
} \over
\o^m(X_1,..., X_m)\cdot  {
|Y(f)|^2 \over | z|^{2d}
}
}
\cdot
{\o^{m+1}\over \si^*\o^{m+1}
}
\right)
\si^*\o^m
\cr
B&=\
\I_Z \log
\left(
{ \o^m(X_1,..., X_m) \cdot  {
|Y(f)|^2 \over | z|^{2d}
}  \over
\o^{m+1}(X_1,..., X_m,Y)
}
\right)
(-{\sqrt{-1}\over 2\pi }\ddb \Phi_\si)\eta_\si
\cr}
$$
Cancelling the common factor $|Y(f)|^2$ in $A$, and using the
simple fact that ${\rm log}\,{\si^*\o^{m+1}\over \o^{m+1}
}=-(m+2)\Phi_\si$, we obtain
$$ A \ = \
\I_Z \ \left\{
\log
\left(
{\si^* \o^m
 \over
\o^m
}
\right)
\ + \
(m+2-d)\Phi_\si
\right\} \cdot
\si^*\o^m
$$
To evaluate $B$, we integrate by parts. Since $Z$ is smooth,
we have
$$
B
=
\I_Z -{\sqrt{-1}\over 2\pi }\ddb\log
\left(
{ \o^m(X_1,..., X_m) \cdot  {
|Y(f)|^2 \over | z|^{2d}
}  \over
\o^{m+1}(X_1,..., X_m,Y)
}
\right)
 \Phi_\si\sum_{i=0}^{m-1}\o^i\o_{\sigma}^{m-1-i}
\eqno(5.7)
$$
Now we clearly have
$$
\ -{\sqrt{-1}\over 2\pi }\pl_X\bar\pl_X\log\left\{
{|Y(f)|^2\over  | z|^{2d} }\right\}  =
d\cdot\o
\eqno(5.8)
$$
since $Y$ is a transversal holomorphic vector field
and $Y(f)$ is a
non-vanishing holomorphic
function on $Z$. Also
$-{\sqrt{-1}\over 2\pi }\pl\bar\pl\log\o^{m+1}
=(m+2)\o$ and
$-{\sqrt{-1}\over 2\pi }\pl_X\bar\pl_X\log\,(\sigma^*\o^m)
=Ric(\sigma^*\o^m)$ by definition.
Thus
$$
B \ = \
\I_Z
\Phi_\si
\Bigg( Ric(\o)
\ + \
d\,\o
\ - \
(m+2)\o
\Bigg)
\cdot
\s_{i=0}^{m-1} \o^i\o_\si^{m-1-i}
$$
This calculation of the $B$ term is very similar to a curvature
calculation in Tian [T1] (see also Lu [Lu2]). Finally,
the remaining contributions from $\log{||\sigma\cdot \Chow(X)||_{\#}^2
\over ||\Chow(X)||_{\#}^2}$ only involve terms of the
form ${||\sigma\cdot \Chow(X)||^2
\over ||\Chow(X)||^2}$, which are known from
Theorem 5
$$
-{D\over V(m+1)}
\log {||\sigma\cdot \Chow(X)||^2
\over ||\Chow(X)||^2}
=F_{\o}^0(\Phi_{\sigma})
=
-{1\over V(m+1)}
\int_Z\Phi_{\sigma}\sum_{i=0}^{m-1}\o^i\o_{\si}^{m-1-i}
$$
Assembling all terms gives Theorem 6. Q.E.D.

\v\v

{\bf \S 6. The Mabuchi energy: General Case}
\v

In this section we establish Theorem 1 for
arbitrary codimension.
Our notation is the one introduced in \S 1.
\v

Let $X \sub \P^N$ be a smooth variety of dimension $n$
and let $Z \sub Gr$ be the corresponding Chow variety, and
$f \in H^0(Gr, O(d))$ be a defining section for $Z$.
Then $Z$ is a singular variety. We say that
the embedding $X\sub \P^N$ is {\it generic} if it satisfies the
following conditions:
\v
1. If $Z_s \sub Z$ is the singular set of $Z$, then $Z_s$
has codimension one.
\v
2.
There exists a subvariety
$Z_{ss} \sub Z_s$ of codimension at least one
(and hence $Z_{ss} \sub Gr$ has codimension at
least two) such that
$Z_s\backslash Z_{ss}$ is a divisor with normal
crossings. In other words,
for all $z \in
Z_s\backslash Z_{ss}$, there exist local coordinates
$(z_0,...,z_m)$, centered at $z$,
 such that $f(z_0,...,z_m) = z_0z_1$.
\v
3. Let $Z_0 = Z\backslash Z_s$ be the set of smooth points
of $Z$. Then the scalar curvature $s(\o)$ is
$L^1(Z_0)$ with respect to the volume form $\o^n$.
\v

\v
The proof consists of two parts. In the first part,
the proof is identical to that for hypersurfaces in the preceding
section, if we view the arguments given there as applying
to the Chow variety $Z$, which is a codimension $1$ subvariety
of the Grassmannian. The result of this first part is to
express $\log\,{||\si\cdot\Chow(X)||_{\#}^2\over
||\Chow(X)||_{\#}^2}$ as the
generalized Mabuchi energy $\nu_{\o_Z}^{\#}(\Phi_{\sigma})$
associated to the singular variety $Z$
(see Lemma 6.1 below). The second part of the proof consists
in identifying the regular part of the generalized Mabuchi
functional $\nu_{\o_Z}^{\#}(\Phi_{\si})$ with the
Mabuchi energy $\nu_{\o_X}(\phi_{\si})$ associated to
the projective variety $X$ itself
(see Lemma 6.2).

\v
Let $\o = \o_Z$ be the restriction of $\o_{Gr}$ to
$Z_0$, the set of smooth points of $Z$. Then
$Ric(\o_Z)$, the Ricci curvature
of $\o_Z$, is a smooth
$(1,1)$ form on $Z_0$. Let $V = V(Z)$ be the volume of $Z$,
and
let $\G_X = \{(x,z) : z \in Gr, x \in X\cap z\}$.
Let $p_1: \G_X: \ra X$ and $p_2: \G_X\ra Z$ be the projection
maps, and let
$Y_s = p_2^{-1}(Z_s) \sub \G_X$ and
$Y_{ss} =  p_2^{-1}(Z_{ss})$.

\v
Since $|s(\o_Z)|$ is
in $L^1(Z_0)$, we can set
$\m(Z) \ = \ {1\over V}\I_{Z_0} s(\o_Z)\o_Z^m$
and define
the {\it Mabuchi functional $\nu_{\o_Z}(\Phi_{\si})$
associated with the regular variety $Z_0$}
as before by
$$
\eqalign{
\n_{\o_Z}(\Phi_\si)
=
{1\over V}
\int_{Z_0} \Bigg\{\log\left({\si^*\o^m_Z\over \o_Z^m}\right)
\si^*\o^m_Z
-\Phi_\si\Big(&  Ric(\o_Z)\s_{i=0}^{m-1}\o_Z^i\si^*\o^{m-1-i}_Z
\cr
&\quad
-
{\m(Z)\over m+1}\s_{i=0}^m
\o_Z^i\si^*\o^{m-i}_Z\Big)\Bigg\}\cr}
\eqno(6.1)
$$
Associated to the variety $Y_s$ is a closed current  $[Y_s]$, supported
on $Y_s$, which
is defined
by the following equation of currents on $\G_X$:
$$-{\sqrt{-1}\over 2\pi} \ddb\log
\left(
{\o_Z^m \wedge \ddb {
|f(z)|^2 \over |P\ell (z)|^{2d}
} \over
\o_Z^{m+1}
}
\right) \ = \ Ric(\o_Z)-(m+2-d)\o_Z -[Y_s]
\eqno(6.2)
$$
\v
The {\it generalized Mabuchi energy $\nu_{\o_Z}^{\#}(\Phi_{\si})$
associated to the singular variety $Z$} can now be defined by
$$
\nu_{\o_Z}^{\#}(\Phi_{\si})
=
\ \n_{\o_Z}(\Phi_\si) \ + \
{1\over V}
\langle [Y_s],\Phi_\si
\s_{i=0}^{m-1}\o_Z^i\si^*\o^{m-1-i}_Z\rangle \
- \ \ {D\over V}\cdot{m\deg(Y_s)\over m+1}\cdot\log
{||\sigma\cdot f||^2
\over ||f||^2}
\eqno(6.3)
$$
where $||\cdot ||$ is the norm introduced in (4.1) and
$\deg(Y_s) = {1\over V}\,\langle[Y_s],\o_Z^{m-1}\rangle $.
\v\v

{\bf Lemma 6.1} {\it Under the above hypotheses, we have
$$
\n_{\o_Z}^\#(\Phi_\si)=
{D(m+2)(d-1)\over V(m+1)}\,{\rm log}\,{||\sigma\cdot f||_{\#}^2
\over ||f||_{\#}^2}
\eqno(6.4)
$$
}
\v\v

{\it Proof.}
Since the map $p_2: \G_X \ra Z$ is birational,
we can pull back integrals over $Z$ to integrals over $\G_Z$
and express
the right side of (6.4) as $A'-B'+T$ where
$$\eqalign{
A'\ = \ &{1\over V}\cdot\I_{\G_X} \log
\left(
{\si^* \o^m \wedge \ddb {
|f(z)|^2 \over |P\ell(\si z)|^2
} \over
\si^*\o^{m+1}
}
\right)
\si^*\o^{m}\cr
B' \ = \ &
{1\over V}\cdot\I_{\G_X} \log
\left(
{\o^m \wedge \ddb {
|f(z)|^2 \over |P\ell( z)|^2
} \over
\o^{m+1}
}
\right)
\o^{m}\cr
}
\eqno(6.5)
$$
and
$$ T \ = \ (m+2-d)\cdot F^o_{\o}(\Phi_\si) \ = \
 -{(m+2-d)\over (m+1)}\cdot{1\over V}\I_{\G_X}
\Phi_\si\cdot\Big(\s_{i=0}^m\o^i\wedge\si^*\o^{m-i}\Big)
$$
where we write $\o$ for $\o_Z$.

\v
Replacing \ $\o^m$ by $ (\o^m - \si^*\o^m) + \sigma^*\o^m$ in (6.5),
we obtain $A'-B'=A-B$ where
\v
$$
\eqalign{
A \ = \ &{1\over V}\I_{\G_X} \log
\left(
{\si^* \o^m \wedge \ddb {
|f(z)|^2 \over |P\ell(\si z)|^2
} \over
 \o^m \wedge \ddb {
|f(z)|^2 \over |P\ell(z)|^2
}
}
\cdot
{\o^{m+1}\over \si^*\o^{m+1}
}
\right)
\si^*\o^m
\cr
& \cr
B \ = \ &{1\over V}\I_{\G_X} \log
\left(
{\o^m \wedge \ddb {
|f(z)|^2 \over |P\ell( z)|^2
} \over
\o^{m+1}
}
\right)
(-{\sqrt{-1}\over 2\pi }\ddb \Phi_\si)
\s_{i=0}^{m-1}
\o^i\wedge\si^*\o^{m-1-i}\
 \cr}
\eqno(6.6)
$$
The same argument as that used in \S 5 again gives
$$
A \ = \
{1\over V}\I_{\G_X} \ \left\{
\log
\left(
{\si^* \o^m
 \over
\o^m
}
\right)
\ - \
(m+2-d)\Phi_\si
\right\} \cdot
\si^*\o^m
\eqno(6.7)$$
In order to evaluate $B$,
we integrate
(6.6) by parts:

$$
B \ = \
{1\over V}\I_{\G_X}\
\Phi_\si\cdot
\Big( Ric(\o)
\ - \
(m+2-d)\o \ - \  [Y_s]
\Big)
\cdot
\s_{i=0}^{m-1}
\o^i\wedge\si^*\o^{m-1-i}
\eqno(6.8)
$$
Assembling (6.6), (6.7) and (6.8) and making use of
$$
m(m+2-d) \ = \ {1\over V}\I_{\G_X} m(Ric(\o)-[Y_s])\o^{m-1} \ = \
\m(Z) - m\deg(Y_s)
$$
we obtain (6.4).
\v
{\it Remarks on $[Y_s]$}
\v
1.
The fact that $[Y_s]$ is supported on $Y_s$ follows from
(5.8).
\v
2. It's not difficult to see that $[Y_s]$ is given by integration
over the smooth points of the variety $Y_s$.
\v
3. The current $[Y_s]$ is also defined by
the equation of currents: $\ddb \log |\nabla F|^2 \ = \ [Y_s]$,
where $\nabla F$ is the gradient of the holomorphic function
which locally defines $Z$.
To see this, we let $\xi$ be the function
$$ \xi(z) \ = \ {\o^m \wedge \ddb {
|f(z)|^2 \over |P\ell( z)|^2
} \over
\o^{m+1}}
$$

Then $\xi$ is a smooth function on $Gr$. Let
$y_0\in \G_X$,
and choose a coordinate system of
$(w_0,...,w_m)$ of a neighborhood of
$p_2(y_0) \in Gr$. Then
$$\xi(y) \ = \ (\xi\circ p_2) (y) \ = \ \s_{i=0}^m
  \  \left| {\pl F\over \pl w_i }
  \right|^2\big(p_2(y) \big)  \ = \
\s_{i=0}^m
  \  \left| f_i(y)\right |^2
$$
in some coordinate neighborhood of $y_0 \in \G_X$, where
$F$ is an analytic function whose divisor is $Z$ and
the $f_i(y)$ are analytic functions whose set of common
zeros is precisely $Y_s$.

\v

\v\v

{\bf Lemma 6.2} {\it Assume that $X\sub \P^N$ is generic. Then
$$ \n_{\o_Z}(\Phi_\si) \ = \ \n_{\o_{FS}}(\phi_\si)
\eqno(6.9)
$$
}
\v

{\it Proof.} As in the proof of Theorem 5,
we differentiate both sides of (6.9)
with respect to $t$ and show that the two sides are equal.
Thus let $c$ be an arbitrary $(N+1)\times(N+1)$
traceless matrix, and let $\si(t)={\rm exp}(ct)\si_0$
for some fixed $\si_0\in G$.
\v
Recall that $\n_{\o_Z}(\Phi_\si)$
is given by (6.1).
Let $\psi_\si = \log\left({\si^*\o^m_Z\over \o_Z^m}\right)$,
and note that
$$ (\si^*\o)^m \cdot{d\over dt}\psi_\si \ = \
(\si^*\o)^m \cdot{m(\si^*\o)^{m-1} \ddb \dot\Phi_\si\over
(\si^*\o)^m \cdot }
\ = \ m(\si^*\o)^{m-1} \ddb \dot\Phi_\si
$$
so
$$
\I_{\G_X} (\si^*\o)^m \cdot{d\over dt}\psi_\si \ = \ 0
$$
Writing $\o_\si = \si^*\o$, we obtain
$$
{d \over dt} \I_{\G_X} \psi_\si (\si^*\o)^m \ = \
\I_{\G_X}\psi_\si m(\si^*\o)^{m-1} \ddb \dot\Phi_\si \ = \
-\I_{\G_X} m\dot\Phi_\si(Ric(\o_\si) - Ric(\o))\o_\si^{m-1}
\eqno(6.10)
$$
\v

On the other hand,
$$ {d \over dt}\I_{\G_X} \s_{i=1}^m
\Phi_\si(\si^*\o)^{m-i}\o^{i-1}Ric(\o)
 \ = \ \I_{\G_X}
\s_{i=1}^m\dot\Phi_\si(\o_\si)^{m-i}\o^{i-1}Ric(\o)
\ + \
$$
$$ \I_{\G_X}  \s_{i=1}^{m-1}\dot \Phi_\si
(m-i)\Big((\o_\si)^{m-i}\o^{i-1}-(\o_\si)^{m-i-1}\o^{i}\Big)
Ric(\o) \ = \
$$
$$\I_{\G_X}
\s_{i=1}^m\dot\Phi_\si(\o_\si)^{m-i}\o^{i-1}Ric(\o) +
\I_X  \s_{i=1}^{m-1}\dot \phi_\si
(m-i)\Big((\o_\si)^{m-i}\o^{i-1}-(\o_\si)^{m-i-1}\o^{i}\Big)
Ric(\o) \ =
$$

$$\I_{\G_X}
\s_{i=1}^m(m-i+1)\dot\Phi_\si(\o_\si)^{m-i}\o^{i-1}Ric(\o)
-
\I_{\G_X}  \s_{i=2}^{m}\dot \Phi_\si
(m-i+1)
(\o_\si)^{m-i}\o^{i-1}
Ric(\o) \ = \
$$
$$\I_{\G_X} m\dot\Phi_\si\o_\si^{m-1}Ric(\o)
\eqno(6.11)
$$
Finally
$$ {d\over dt}{\m(Z)\over m+1}\cdot\
\I_{\G_X}\s_{i=0}^m\o^i\o_\si^{m-i} \ = \
\I_{\G_X} \dot\Phi_\si \m(Z)\o_\si^m
\eqno(6.12)
$$
Combining (6.10), (6.11), and (6.12) we obtain
$$
{d\over dt} \n_{\o_Z}(\Phi_\si) \ = \
-{1\over V}\cdot \I_{\G_Z} \dot\Phi_\si
(m Ric(\o_\si) - \m(Z)\o_\si)\o_\si^{m-1}
\eqno(6.13)
$$

\v
This is the derivative of the Mabuchi energy on $Z$.
In order to establish Lemma 6.2 we must show
it is also the derivative of the
Mabuchi energy on $X$:
 Replacing $X$ and $Z$ by
$\si_0(X)$ and $\si_0(Z)$, we see that we may assume that
$\si_0$ is the identity matrix. We shall do that from
now on, and we shall write $\dot\Phi_{\sigma_0}(c) =
\dot\Phi(c)$ and
$\dot\phi_{\sigma_0}(c) =
\dot\phi(c)$.
\v
Thus, we let
$$ M_Z \ = \
{\sqrt{-1}\over 2\pi } Ric(\o_{Gr}^m) \ -
{\m(Z)\over m}\o_{Gr}
$$
and
$$ M_X \ = \
{\sqrt{-1}\over 2\pi } Ric(\o^n) \ -
{\m(X)\over n}\o
$$
We claim that:
$$ {1\over vol(X)}\I_X M_X\ \dot\phi(c)\  n\o^{n-1}  \ = \
{1\over vol(Z)}\I_Z M_Z\ \dot\Phi(c)\  m\o_{Gr}^{m-1}
\eqno(6.14)
$$
Theorem 1 then will follow from (6.13) and (6.14).
\v
 Let
$\Gamma_X  = \{(x,z) \in \Gamma : x \in X\}$. Let $p_i = \pi_i|_X$.
Then $p_1: \Gamma_X \ra X$ and $p_2: \Gamma_X \ra Z$,
and we have the following double fibration
$$
\matrix{ &\Gamma_X&   \cr
{}^{p_1}\swarrow &    &\searrow{}^{p_2}  \cr
X&                     & Z\cr}
\eqno(6.15)
$$
and that currents on $X$ and $Z$ can then be compared by
pull-backs and push-forths through this fibration.
Note that
$p_2$ is a birational map. This means that there is a Zariski open
subset $Z_0 \sub Z$ such that $p_2: p_2^{-1}(Z_0) \ra Z_0$ is bijective.
Thus
$$
\I_Z M_Z\ \dot\Phi(c)\  m\o_{Gr}^{m-1} \ =
 \
\I_{\G_X}
p_2^* M_Z\cdot  p_2^*\left(\dot\Phi(c)\  m\o_{Gr}^{m-1}\right) \ =
$$
$$\I_{\G_X}
p_1^* M_X\cdot  p_2^*\left(\dot\Phi(c)\  m\o_{Gr}^{m-1}\right) \ + \
\I_{\G_X}
(p_2^* M_Z-p_1^*M_X)\cdot  p_2^*\left(\dot\Phi(c)\  m\o_{Gr}^{m-1}\right)
$$
\v\v
But according to (4.6),
$$ {p_1}_*p_2^*(m\dot\Phi(c)
\sigma ^*\o_{Gr}^{m-1})
 \ = \
n\dot\phi(c)\o^{n-1}
$$
and so we get
$$
\I_{\G_X}
p_1^* M_X\cdot  p_2^*\left(\dot\Phi(c)\  m\o_{Gr}^{m-1}\right) \ = \
\I_X M_X\ \dot\phi(c)\  n\o^{n-1}
$$
Thus (6.14) will follow from:
$$
\I_{\G_X}
M(X,Z)\cdot  p_2^*\left(\dot\Phi(c)\  m\o_{Gr}^{m-1}\right)
\ = \
\I_{X}{p_1}_*\left[
M(X,Z)\cdot  p_2^*\left(\dot\Phi(c)\  m\o_{Gr}^{m-1}\right)
\right] \ =
\ 0
\eqno(6.16)
$$
where
$$ M(X,Z) \ = \ (p_2^* M_Z-p_1^*M_X)
$$
Now $M(X,Z)$ can be made explicit as follows: Define a
function $A_X: \G_X \ra \C$ by the formula
$$  p_2^*\o_{Gr}^m \ = \ A_X\cdot p_2^*\o_{Gr}^{m-n}p_1^*\o^n
$$
Then $M(X,Z) \ = \ {\sqrt{-1}\over 2\pi }R(h_X)$ where $h_X$ is the curvature
of the
metric on the relative canonical bundle
$\cL_X = K_{\G_X}\otimes p_1^*(K_X)^{-1}$ given by the
formula
$$  h_X \ = \ A_X \cdot p_2^*\o_{Gr}^{m-n}
$$
\v

We provide now a proof of the key equation (6.16).
The main idea is that, by making use of a first jet extension
(see the map $\iota_X(x)=(x,\zeta)$ defined below),
the double fibration (6.15) can be imbedded in another
double fibration, which is independent of $X$
and is more symmetric.
More precisely, we shall define manifolds
$\Sigma$ and  $\Gamma'$,  a fibration $P_1:\Sigma \ra \Gamma'$,
a metric
$h$ on the relative canonical bundle
$\cL = K_\Sigma \otimes P_1^*(K_{\G'})^{-1}$ and
embeddings $\iota:\Gamma_X \hookrightarrow \Sigma$,
$\iota_X: X \hookrightarrow \Gamma'$ which will have the following
properties: The projection $p_1: \Gamma_X \ra X$ is the
restriction of $P_1$, metric $h_X$ is the restriction
of the $h$ to $\cL_X$,
$$
\matrix{\Gamma_X & \ra & \Sigma\cr
\downarrow &  & \downarrow \cr
X & \ra & \Gamma'\cr}
\eqno(6.17)
$$
The reason
for doing this is that $\Sigma$ and $\Gamma'$ will have a lot
of symmetries, which will help us compute the fiber integral of $M$.
\v
Here are the definitions:
\v
$Gr = Gr(N-n-1, \P^N)$

$Gr' = Gr(n, \P^N)$

$\Sigma \ = \ \{(x,\z,z): x \in \P^N, \z \in Gr', z \in Gr\ ,\
x \in \z \cap z \}$

$\Gamma = \{ (x,z): x \in \P^N, z \in Gr, x \in z \}$

$\Gamma' = \{ (x,\z): x \in \P^N, \z \in Gr', x \in \z \}$

$P'(x,\z,z) = (x, \z) $

$P(x,\z,z) = (x, z) $

$\Gamma_X = \{ (x,z): x \in X, z \in Gr, x \in z \}$

$\iota_X(x) = (x, \z)$ where $\z$ is the unique hyperplane
of dimension $n$ tangent to $X$ at $x$.
\v\v

Note that $G=SL(N+1)$ and $U = U(N+1)$ act on $\Sigma, \Gamma,
\Gamma', Gr, Gr'$, that $U$ leaves the metrics
$\o, \o_{Gr}, \o_{Gr}'$ invariant, and that the original double
fibration (6.15) has now been extended to the following double
fibration
$$
\matrix{ &\Sigma_X&   \cr
{}^{P'}\swarrow &    &\searrow{}^{P}  \cr
\Gamma' &                     & \Gamma\cr}
\eqno(6.18)
$$
\v
Define a function  $A: \Sigma \ra \C$
as follows:

$$ A(x, \z, z) \ = \ {\o_{Gr}^m (X_1,...,X_n, Y_1,...,Y_{m-n})
\over \o^n(X_1,...,X_n)\o_{Gr}^{m-n}(Y_1,...,Y_{m-n})}
$$

where $a=(x,\z,z)$,  $Y_1,...,Y_{m-n}$ is a basis for $T_a(F_a)$,
and
$X_1,...,X_n$ are tangent vectors in $T_a(\Sigma)$ whose projections
to $\P^N$ form a basis of $T_x(\z)$. This function
is clearly invariant under the action of $U$.
\v
Define a metric $h$ on $\cL$ by the formula:

$$  h \  = \ A \cdot P_2^*\o_{Gr}^{m-n}
$$
where $P_2: \Sigma \ra Gr$ is the projection map.
Let
$$M = {\sqrt{-1}\over 2\pi} R(h)$$

\v
Then (6.16) is equivalent to
$$
\I_{X} \iota_X^*{P_1}_*\left[
M\cdot  P_2^*\left(\dot\Phi(c)\  m\o_{Gr}^{m-1}\right)
\right] \ =
\ 0
\eqno(6.19)
$$
In fact, we shall prove that
$$
\iota_X^*{P_1}_*\left[
M\cdot  P_2^*\left(\dot\Phi(c)\  m\o_{Gr}^{m-1}\right)
\right] \ =
\ 0
\eqno(6.20)
$$
Let $T(\G')$ be the tangent bundle of $\G'$ and let
$S \sub T(\G')$ be the subundle defined by
$$  S \ = \ Ker(T(\G') \ra T(Gr'))
$$
Thus $S$ is a vector bundle on $\G'$ of rank $n$ and
for $(x,\z) \in \G'$,
the fiber $S_{(x,\z)}$ is the tangent space of $\z$
at the point $x$.
Relation (6.20) then follows from the following:

$$
{P_1}_*\left[
M\cdot  P_2^*\left(\dot\Phi(c)\  m\o_{Gr}^{m-1}\right)
\right]\Big|_S \ =
\ 0
\eqno(6.21)
$$
Thus  we must show,
for every $(x,\z) \in \G'$, the
$(n,n)$ form
${P_1}_*\left[
M\cdot  P_2^*\left(\dot\Phi(c)\  m\o_{Gr}^{m-1}\right)
\right]$
evaluated at a generator of
$\Lambda^n T_x(\z)\otimes \Lambda^n T'_x(\z)$, is zero
(here $T$ is the holomorphic tangent space and $T'$
the anti-holomorphic tangent space).
\v
Thus we fix $(x,\z) \in \G'$. Recall
that $x\in \P^N$, that is, $x\sub \C^{N+1}$. Also,
we have $x \in \z$ where $\z$ is a plane in
$\P^N$ of dimension $n$, that is $\z \sub \C^{N+1}$ is a
vector space of dimension $n+1$. The tangent space
of $\z$ at the point $x$ is canonically isomorphic
to $x^\perp \cap \z = \z' \sub \C^{N+1}$. Let
$\z_1,...,\z_n$ be an orthonormal basis of
$\z'$. Let

$$ B(c,x;\z_1,...,\z_n) \ = \
{P_1}_*\left[
M\cdot  P_2^*\left(\dot\Phi(c)\  m\o_{Gr}^{m-1}\right)
\right]\Big|(\z_1\wedge \bar \z_1, ... , \z_n\wedge \bar \z_n)
$$

We want to show that $B(c,x;\z_1,...,\z_n) = 0$.
\v
Recall that

$$ \dot\Phi(c) \ = \ tr((c+c^*)ZZ^*)
$$
where, abusing notation, $Z$ is an $(N+1) \times (N-n)$
matrix whose columns form an orthonormal basis of
the vector space $Z \sub \C^{N+1}$.
\v
Define an $(N+1) \times (N+1)$ matrix
$$B(x,\z)\ =  \ B(x;\z_1,...,\z_n) \ = \
{P_1}_*\left[
M\cdot  P_2^*\left(ZZ^*\  m\o_{Gr}^{m-1}\right)
\right]\Big|(\z_1\wedge \bar \z_1, ... , \z_n\wedge \bar \z_n)
\eqno(6.22)
$$
Then $ B(c,x;\z_1,...,\z_n) = tr\big((c+c^*)B(x,\z)\big)$
so it suffices to show that $B(x,\z) = 0$ for all $(x,\z) \in\G'$.
\v
Note that $B$ has the following properties:

$$  B(u(x,\z)) \ =  \ u B(x,\z)u^*, \ \ B(x,\z u_1) \ = \ B(x,\z)
$$
for all $u \in U(N+1)$ and all $u_1 \in U(n)$.
This implies that

$$  B(x,\z) \ = \ k\cdot \z \z^*
$$
for some constant $k$.
We claim that $k=0$. To see this, we take the trace of (6.22)
and use the fact that $tr(\z\z^*) = tr(\z^*\z) =  (n+1)$
since $\z\z^*$ is the identity matrix, and that likewise,
$tr(ZZ^*) = tr(Z^*Z) = (N-n)$:

$$  k\cdot (n+1) \pi^*\o^n_{FS} \ = \
{P_1}_*\left[
M\cdot  P_2^*\left((N-n)\  m\o_{Gr}^{m-1}\right)
\right]\Big|(\z_1\wedge \bar \z_1, ... , \z_n\wedge \bar \z_n)
\eqno(6.23)
$$
where $\pi: \G' \ra \P^N$ is the projection map.
Applying $\iota_X^*$ to both sides and integrating
over $X$ we get

$$ k\cdot (n+1) \cdot \I_X \o^n_{FS} \ = \ (N-n)\I_{\G_X}
(p_2^*M_Z - p_1^*M_X) p_2^*\o_{Gr}^{m-1}
$$
But
$$ \I_{\G_X} p_2^*M_Z\cdot  p_2^*\o_{Gr}^{m-1} \ = \
\I_Z M_Z \cdot \o_{Gr}^{m-1} \ = 0
$$
Also, making use of (4.7) we have
$$ \I_{\G_X} p_1^*M_X \cdot p_2^*\o_{Gr}^{m-1} \ = \
\I_{X}  M_X \cdot {p_1}_*p_2^*\o_{Gr}^{m-1} \ = \
D_1n \I_{X}  M_X \cdot \o^{n-1}_{FS} \ =  0
$$
Thus $k\cdot (n+1) = 0 - 0$ so $k=0$. Thus $B(x,\z) = 0$ and
the proof of Lemma 6.2,
and hence of Theorem 1 is complete. Q.E.D.

\v\v

{\it Remark.} Although we do not require it in the preceding proof,
it may be useful to note that the key function $A(x,\zeta,z)$
introduced can be described by a simple explicit formula.
Let $x^{\perp}\sub \C^{N+1}$ be the orthogonal
complement of $x$, and let $\pi_x:\C^{N+1}\rightarrow x^{\perp}$
be the orthogonal projection. Then
$\pi_x(\zeta)\sub x^{\perp}$
has dimension $n$ and $\pi_x(z)\sub x^{\perp}$ has dimension
$N-n-1$. These two spaces will generically span a
subspace $\pi_x(\zeta+z)$ of codimension $1$ inside $x^{\perp}$.
Then $A(x,\zeta,z)$ is the length of the Pl\"ucker vector
of $\pi_x(\zeta+z)$ with respect to a basis which is the
union of an orthonormal basis for $\pi_x(z)$ and $\pi_x(\zeta)$
$$
A(x,\zeta,z)
=
{|Pl(\pi_x(\zeta+z))|\over
|Pl(\pi_x(\zeta))|\cdot |Pl(\pi_x(z))|}
\eqno(6.24)
$$
Another way of writing this is as follows. Fix $x$ and let $H=x^{\perp}$.
View $A$ as a function on
$Gr(n,H)\times Gr(N-n-1,H)$. Let
$$
\theta:
Gr(n,H)\times Gr(N-n-1,H)
\rightarrow Gr(N-1,H)
$$
be the map which sends $(\zeta_0,z_0)$ to the space spanned by
$\zeta_0$ and $z_0$. Then
$$
-\pl\bar\pl{\rm log}\,A=
\theta^*\o_{N-1}-\o_n-\o_{N-n-1}
$$
where the $\o_k$ are just Fubini-Study metrics. In other words,
$R(h)=\theta^*\o_{N-1}-\o_n$.

\v\v
{\it Remark.} Generalized Mabuchi energies for singular varieties
emerge naturally from the above proof, and may be worth
investigating in their own right. Related extensions for
the Futaki invariant were studied in Ding and Tian [DT].

\v\v
{\it Remark.} There are many ways of expressing the current $[Y_s]$.
It is intriguing that the term it produces can be formally
viewed as an $F^0$ functional for $Y_s$.

\v\v
{\it Remark.} The genericity assumption on $X$ is not
really restrictive. It should suffice for the study of asymptotic
stability, when the variety $X$ is imbedded into projective space
by the bases of the antipluricanonical bundle $K_X^{-p}$
for $p$ large.

\v\v
{\it Remark.} The seminorm $||\cdot||_{\#}$
is degenerate. Its main property
is that for any norm $||\cdot||_B$ on $H^0(Gr,O(d))$ with $||f||>0$ for $f\not=0$,
there exists a constant $C>0$ so that
$$
||f||_{\#}\leq C ||f||_B
\eqno(6.25)
$$
for all $f\in H^0(Gr,O(d))$. This follows from the continuity of $||\cdot||_{\#}$,
which requires a somewhat technical argument. For our purposes,
it suffices to observe that the second term in the expression (1.1)
for $||\cdot||_{\#}$ is continuous, since the only possible divergences of the
integrand is logarithmic. As for the first term, we can estimate it as follows
$$
\I_Z{\rm log}\,\left(\o_{Gr}^m\wedge \ddb
\left({|f(z)|^2\over |P\ell(z)|^{2d} }\right)
\right)\wedge \o_{Gr}^m
\leq C\,{\rm sup}_{Gr}{\rm log}\,|\nabla f|^2
\I_Z\o_{Gr}^m
\eqno(6.26)
$$
Since the logarithmic is an increasing function
and ${\rm sup}_{Gr}{\rm log}\,|\nabla f|^2={\rm log}\,{\rm sup}_{Gr}|\nabla f|^2$,
the right hand side is bounded on the unit ball with respect to $||\cdot||_B$
in $H^0(Gr,O(d))$, and our claim follows.
\v
{\it Remark.} As mentioned in \S 5, the Deligne Pairing
is related to the concept of CM stability, as defined by
Tian [T2]:
\v
Let $G = SL(N+1,\C)$. Let $\pi:\cX \ra B $ be a $G$
equivariant holomorphic fibration between smooth
varieties, equivariantly embedded in $B\times \P^N$.
Tian constructs a $G$ equivariant line bundle
$L_B$ over $B$ and a metric $||\cdot ||_Q$ on $L_B$
with the following property: Let $b\in B$ and let
$X = \pi^{-1}(b)$.

$$
\n_\o(\phi_\sigma)
= \
C\cdot\,{\rm log}\,{||\sigma\cdot b ||_{Q}^2
\over ||b||_{Q}^2}
\eqno(6.27)
$$
for some positive constant $C$.
\v

Now we let $\cK$ be the relative canonical bundle of
$\cX$ over $B$ and let $\cL$ be the $O(1)$ bundle
on $\cX$. Let
$$\cM \ =\  <\cK,\cL,\cL,...,\cL>(\cX/B)
$$
be the Deligne pairing of the bundle $\cK$ and $n$
copies of the bundle $\cL$. Then
$\cM$ is a $G$ equivariant line bundle
$L_B$ over $B$ and comes equipped with
the Deligne metric $||\cdot ||_D$.
Let $b\in B$ and let
$X = \pi^{-1}(b)$. Then we can show that

$$
\n_\o(\phi_\sigma)
= \
C\cdot\,{\rm log}\,{||\sigma\cdot b ||_{D}^2
\over ||b||_{D}^2}
\eqno(6.28)
$$
It is natural to compare (6.27) and (6.28) and
investigate their relationship. Also,
 the role of the pairing
$<\cK,\cK,...,\cK,\cL,\cL,...,\cL>(\cX/B)$
($p$ copies of $\cK$ and $q$ copies of $\cL$
with $p+q=n+1$) needs to be clarified.
These questions  will be addressed in an
upcoming paper.

\v\v\v

{\bf ACKNOWLEDGEMENTS}

\v

The authors would like to thank  Gang Tian
for many very helpful suggestions, especially
on singularities, and  Sean Paul
for calling their attention to important examples of
unstable manifolds. They would also like to thank
J. Song, Mu-Tao Wang, Ben Weinkove, and all
the members of the Columbia geometric analysis group
for many stimulating exchanges.

\vfill \eject
\bigskip

\centerline{\bf REFERENCES}

\v\v

[BM] Bando, S. and T. Mabuchi,
``Uniqueness of Einstein-K\"ahler metrics modulo
connected group actions",
in {\it Algebraic Geometry, Sendai 1985},
Adv. Stud. in Pure Math. {\bf 10} (1987)
11-40, Kinokuniya, Tokyo and North-Holland,
Amsterdam.

\v

[C1] Chen, X.,
``On the lower bound of the Mabuchi energy
and its application",
Int. Math. Res. Notices {\bf 12} (2000) 607-623.

\v
[C2] Chen, X.,
``The space of K\"ahler metrics", J. Differential Geom.
{\bf 56} (2000) 189-234.

\v

[De] Deligne, P.,
``Le determinant de la cohomologie",
Contemporary Math. {\bf 67} (1987) 93-177.
\v

[DT] Ding, W. and G. Tian,
``K\"ahler-Einstein metrics and the generalized Futaki
invariants", Inventiones Math. {\bf 110} (1992) 315-335.

\v

[D1] Donaldson, S.,
``Anti self-dual Yang-Mills connections over complex
algebraic surfaces and stable vector bundles",
Proc. London Math. Soc. {\bf 50} (1985) 1-26.

\v
[D2] Donaldson, S.,
``Scalar curvature and projective imbeddings I",
2001 Preprint.

\v
[D3] Donaldson, S.,
``Infinite determinants, stable bundles, and curvature",
Duke Math. J. {\bf 54} (1987) 231-247.

\v
[KN] Kempf, G. and L. Ness,
``The length of vectors in representation spaces",
{\it Algebraic Geometry} (Proc. 1978
Summer Meeting, Copenhagen), Lecture Notes in Math.
{\bf 732} (1979) 233-243, Springer-Verlag.
\v

[Lu1] Lu, Z.,
``On the lower order terms of the asymptotic expansion
of Tian-Yau-Zelditch", Amer. J. Math. {\bf 122} (2000)
235-273.

\v

[Lu2] Lu, Z.,
``K-energy and K-stability for hypersurfaces",
math.DG/0108009

\v

[Luo] Luo, H.,
``Geometric criterion for Gieseker-Mumford stability
of polarized manifolds",
J. Diff. Geom. {\bf 49} (1998) 577-599.

\v

[M] Mabuchi, T.,
``K-energy maps integrating Futaki invariants",
Tohoku Math. J. {\bf 38} (1986) 245-257.

\v

[Mu] Mumford, D.,
``Stability of projective varieties",
L'Enseignement Mathematique,
{\bf 23} (1977) 39-110.

\v

[P] Paul, S.
``Geometric analysis of Chow Mumford stability",
Princeton Ph.D. Thesis (2000).

\v
[SY] Siu, Y.T. and S.T. Yau,
``Complete K\"ahler manifolds with nonpositive curvature
of faster than quadratic decay",
Ann. of Math. {\bf 105} (1977) 225-264.

\v

[T1] Tian, G.,
``The K-energy on hypersurfaces and stability",
Comm. Anal. Geometry {\bf 2} (1994) 239-265.

\v

[T2] Tian, G.,
``K\"ahler-Einstein metrics with positive scalar
curvature", Inventiones Math. {\bf 130} (1997) 1-37.

\v
[T3] Tian, G.,
``{\it Canonical Metrics in K\"ahler Geometry}",
Birkh\"auser, Basel, 2000.

\v
[T4] Tian, G.,
``Bott-Chern forms and geometric stability",
Discrete Contin. Dynam. Systems {\bf 6} (2000)
211--220.

\v

[UY] Uhlenbeck, K. and S.T. Yau,
``On the existence of Hermitian Yang-Mills connections
on stable vector bundles", Commun. Pure Appl. Math.
{\bf 39} (1986) 257-293.

\v

[W] Wang, X.,
``Balance point and stability of vector bundle over
projective manifolds", Brandeis 2001 preprint.

\v
[Y1] Yau, S.T.,
``On the Ricci curvature of a compact K\"ahler manifold
and the complex Monge-Ampere equation I",
Comm. Pure Appl. Math. {\bf 31} (1978) 339-411.

\v

[Y2] Yau, S.T.,
``Open Problems in Geometry",
Proc. Symposia Pure Math. {\bf 54} (1993) 1-28.

\v
[Y3] Yau, S.T.,
``Review of K\"ahler-Einstein metrics in Algebraic Geometry",
Israel Math. Conf. Proceedings {\bf 9} (1996)
433-443.

\v
[Y4] Yau, S.T.,
``Nonlinear analysis in geometry",
Enseign. Math. {\bf 33} (1987) 109-158.

\v

[Z] Zhang, S.,
``Heights and reductions of semi-stable varieties",
Compositio Math. {\bf 104} (1996) 77-105.

\end